\documentclass{article}

\usepackage{arxiv}

\usepackage{hyperref} 
\usepackage{url} 
\usepackage{booktabs} 
\usepackage{amsfonts} 
\usepackage{nicefrac} 
\usepackage{microtype} 
\usepackage{lipsum}

\usepackage{graphicx}
\usepackage{amsmath}

\usepackage{xcolor}

\usepackage{ulem}

\usepackage[english]{babel}

\title{\textcolor{black}{High-Order Block Toeplitz Inner-Bordering method for solving the Gelfand-Levitan-Marchenko equation}}

\author{
 Sergey Medvedev$^{1,2,*}$, Irina Vaseva$^{1,2}$, Mikhail Fedoruk$^{2,1}$\\
$^{1}$ Federal Research Center for Information and Computational Technologies,\\ Novosibirsk
630090, Russia,\\
$^{2}$ Novosibirsk State University, Novosibirsk 630090, Russia,\\
* Corresponding author: medvedev@ict.nsc.ru
 }

\begin{document}
\maketitle

\begin{abstract}
\textcolor{black}{We propose a high precision algorithm for solving the  Gelfand-Levitan-Marchenko equation. The algorithm is based on the block version of the Toeplitz Inner-Bordering algorithm of Levinson's type. To approximate integrals, we use the high-precision one-sided and two-sided Gregory quadrature formulas. Also we use the Woodbury formula to construct a computational algorithm. This makes it possible to use the almost Toeplitz structure of the matrices for the fast calculations.} 
\end{abstract}

\keywords{Gelfand-Levitan-Marchenko equations \and Toeplitz Inner-Bordering method \and Zakharov-Shabat system \and Inverse spectral problem \and \textcolor{black}{Gregory quadrature formulas} \and Woodbury formula}

\section{Introduction}

The inverse scattering transform for the nonlinear Schrödinger equation (NLSE) allows one to integrate this 
equation~\cite{ZakharovShabat1972, ablowitz1981solitons,lamb1980elements}. This method, also known as nonlinear Fourier transform (NFT), consists of two parts: the direct and the inverse ones. The first part is the direct spectral Zakharov-Shabat (ZS) problem, in which the scattering data is determined. Scattering data, also referred as a nonlinear spectrum, consists of continuous and discrete spectra. Discrete spectrum exists only in the case of anomalous dispersion and describes the solitonic solutions of the NLSE. The inverse NFT consists of restoring the NLSE solution from its nonlinear spectrum.

 NFT has recently received a lot of attention in areas where NLSE is used to describe various types of optical signals. In particular, it is used in telecommunications applications, where a new method has been proposed to compensate for the effects on a signal during its propagation in optical fiber~\cite{Yousefi2014III, Le2014, Gui2017a, Wahls2017, Gui2018, Civelli2019}. The method is also used to describe and analyze various physical phenomena~\cite{Gelash2019PRL, Mullyadzhanov2019}. In addition, the problem of scattering on Bragg gratings, which are the basis for optical filters in high-speed fiber optic data lines, is reduced to the ZS system \cite{kashyap1999fibre, podivilov2006exactly}. The linear Schrödinger equation for two-level quantum systems with a time-varying Hamiltonian takes the form of a ZS system \cite{akulin2005coherent,carmel2000geometrical}.

The idea of using solitons for data transmission in fiber optic lines first appeared in the work ~\cite{hasegawa1973transmission}. Since this work, NLSE and its modifications have been intensively studied in connection with fiber optic telecommunication systems\cite{hasegawa2003optical,mollenauer2006solitons,turitsyn2003physics}. Later, the idea was put forward to use multisoliton pulses in fiber optic data lines, when the information is modulated and reconstructed in the so-called nonlinear Fourier space \cite{hasegawa1993eigenvalue,hari2016multieigenvalue}.

In this regard, in recent years, interest in numerical methods for solving direct and inverse problems has increased. This paper focuses on the numerical methods for the inverse NFT, i.e. a signal reconstruction from a scattering data. There are several approaches to the inverse NFT. The common ones are based on the solution of the Riemann–Hilbert problem \cite{wahls2016fast,kamalian2020full} or the Gelfand–Levitian–Marchenko equations (GLME). In \cite{wahls2016fast} it is shown that methods based on the numerical solution of the GLME are more accurate in comparison with the method based on factorization for the Riemann-Hilbert problem. Among the GLME based methods we can mention: a method based on the transition to a system of partial differential equations \cite{Xiao2002}, the Toeplitz inner bordering method (TIB)  \cite{belai2006finite,belai2007efficient,frumin2015efficient} and its generalized block version (GTIB) \cite{MedvedevGTIB}, the integral layer peeling (ILP) method \cite{feced1999efficient,Rosenthal2003} and an algorithm using polynomial interpolation to approximate the GLME kernels \cite{Ahmad1998}. 

In the case of a purely continuous spectrum, time-reversed forward NFT algorithms can be used to implement the inverse NFT \cite{wahls2016fast,Yousefi2020}. In this case, the inverse NFT can be considered as a counterpart to the direct NFT, similar to the conventional Fourier transform. For a purely discrete spectrum, the most effective method is the Darboux transform \cite{aref2016control}. The Darboux transform can also be used in the case of a combination of continuous and discrete spectra \cite{aref2018modulation}, but it requires the solution of the GLME for the continuous part of the spectral data.
An overview of methods for solving the direct and inverse NFT problem can be found, for example, in \cite{turitsyn2017nonlinear}.
Comparison of the efficiency of numerical methods for Bragg lattices is given in \cite{buryak2009comparison}.

The inverse NFT methods mentioned above do not achieve higher than second order accuracy. We propose an approach allowing one to increase the precision of the inverse NFT up to sixth or seventh order. The approach is based on the numerical solution of the GLME using the second order GTIB method \cite{MedvedevGTIB}. To approximate integrals in the GLME, we use the high-precision one-sided and two-sided Gregory quadrature formulas. Also we use the Woodbury formula to construct a computational algorithm. This makes it possible to use the almost Toeplitz structure of the matrices for a fast calculations. The numerical comparison is made with the second order TIB method  \cite{frumin2015efficient}, which has proven to be one of the most efficient methods applied for continuous and discrete spectra.

\section{ Gelfand-Levitan-Marchenko equation} 

The direct NFT for the complex-valued potential (signal) $q(t)$ is computed as a solution of the ZS system \cite{ZakharovShabat1972}
\begin{equation}\label{ZS}
\psi_{1t}+i\zeta \psi_1=q(t)\psi_2,\quad \psi_{2t}-i\zeta\psi_2=\mp q^*(t)\psi_1, \end{equation}
where $t\in\mathbb{R}$, $\psi(t) = [\psi_1,\psi_2]^T \in\mathbb{C}^2$ is a wave function, $\zeta=\xi+i\eta\in\mathbb{C}^+$ is a spectral parameter.
Here the upper sign (minus) corresponds to the case of anomalous dispersion, the lower sign (plus) -- normal dispersion. 
The problem is considered under the assumption that $q(t)$ decays at least exponentially
for $t \rightarrow \pm\infty$.

Solving the ZS system (\ref{ZS}) we can find the left 
NFT spectrum of the signal $q(t)$ \cite{turitsyn2017nonlinear}:
\begin{equation}\label{SpectralData}
\Sigma_l=\left\{l(\xi),\left[\zeta_n,l_n\right]_{n=1}^N\right\}.
\end{equation}
The NFT spectrum consists of the continuous and discrete spectra. The continuous spectrum is defined by the left $l(\xi)$ 
reflection coefficient
with respect to the real-valued spectral parameter $\xi$. The discrete spectrum is presented by
$N$ eigenvalues $\zeta_n$ of the ZS system, corresponding to the solitons of the signal $q(t)$, and
left $l_n$ 
norming constants (also referred as phase coefficients). 
Note that in the case of normal dispersion the spectrum cannot have the discrete part. 

To recover the potential $q(t)$ from the spectral data (\ref{SpectralData}) one need to solve
the inverse problem for the ZS system (\ref{ZS}). It can be
 reduced to solving the left system of integral equations \cite{lamb1980elements}
\begin{equation}\label{A1}
\begin{array}{l}
A_1^*(t,s)+\int\limits_{-\infty}^t\,A_2(t,t')\,\Omega_l(t'+s)\,dt'=0,\quad t\geq s,\\
\mp A_2^*(t,s)+\Omega_l(t+s)+\int\limits_{-\infty}^t\,A_1(t,t')\,\Omega_l(t'+s)\,dt'=0.
\end{array}
\end{equation}
where the kernel $\Omega_l(t)$ is defined for all real $t$ by
\begin{equation}\label{Omega_l}
\Omega_l(t)=\frac{1}{2\pi}\int\limits_{-\infty}^\infty l(\xi)e^{-i\xi t}\,d\xi -i\sum\limits_{n=1}^N\,l_{n}e^{-i \zeta_n t}.
\end{equation}
The equations (\ref{A1}) 
are referred as the left 
Gelfand-Levitan-Marchenko equations (GLME). Further we omit the subscript $l$ in the notation of the kernel $\Omega_l$.

After solving the system (\ref{A1}), 
the potential $q(t)$ is restored by the formulas
\begin{equation}\label{q=2A2}
q(t)=-2A_2^*(t,t).
\end{equation}




Let's make the change of variables: $k=t-s$, $n=t$,
\begin{equation}\label{knts}
X_1(k,n)=A_1(t,s),\quad X_2(k,n)=\mp A_2^*(t,s).
\end{equation}
Taking into account the condition $t\geq s$, we obtain that the functions $X_{1,2}(k,t)$ 
are defined on the half-space $k\geq 0$.
Then the equations (\ref{A1}) take the form
\begin{equation}\label{X1}
\begin{array}{l}
X_1(k,t)\mp\int\limits_0^{\infty}\,\Omega^*(2t-k-p)X_2(p,t)\,dp=0,\\
X_2(k,t)+\int\limits_0^{\infty}\,\Omega(2t-k-p)X_1(p,t)\,dp=-\Omega(2t-k).
\end{array}
\end{equation}
For fixed $t$ the integral operators are Hankel, since they depend only on the sum $k+p$.

Let's make the first approximation. It consists in replacing the infinite domain of integration with a finite one $ [0, P] $ for a sufficiently large value $P$. We obtain integral equations with finite limits ($0\leq k\leq P$)
\begin{equation}\label{X1P}
\begin{array}{l}
X_1(k,t)\mp\int\limits_0^{P}\,\Omega^*(2t-k-p)X_2(p,t)\,dp=0,\\
X_2(k,t)+\int\limits_0^{P}\,\Omega(2t-k-p)X_1(p,t)\,dp=-\Omega(2t-k).
\end{array}
\end{equation}

{\textcolor{black}{
We will use the Hankel structure of integral operators to construct a fast algorithm to solve the system (\ref{X1P}). Firstly, since after transformations the potential is determined as $q(t)=\pm 2X_2(0,t)$, the point $k=0$ must be included in the computational grid.
Secondly, we apply a grid with a constant step size $h$ and divide the interval $[0,P]$ on $M$ subintervals of the length $h$ and as result $P=M h$. This allows us to approximize the kernel $\Omega(2t-k-p)$ of the integral operators  by the Hankel matrixes $H$ and $H^*$
\begin{equation}\label{matrix_H}
H =  \begin{bmatrix}
\omega_0& \omega_1               & \omega_2    &\cdots                  &\omega_M\\
\omega_1& \omega_2               &             &\rotatebox{70}{$\ddots$}&        \\
\omega_2&                        &\rotatebox{70}{$\ddots$}&             &\vdots  \\ 
\vdots  &\rotatebox{70}{$\ddots$}&             &                        &\omega_{2M-1}\\
\omega_M&                        & \cdots      &  \omega_{2M-1}         &\omega_{2M}
\end{bmatrix},\quad \omega_k=\Omega(2t-k h).
\end{equation}
}}
{\textcolor{black}{
Thirdly, we use the Gregory's formula for numerical integration \cite{hamming1962numerical, phillips1972gregory}. These quadrature formulas contain a minimum number of weight coefficients different from unit. And these coefficients are located symmetrically on the edges of the integration interval. For our calculations, it is convenient to present weight coefficients from the Gregory formula in the form of diagonal matrices $W_n$ with size $M+1$
\begin{equation}
W_1=\mbox{diag}\left[\frac{1}{2},1,\ldots,1,\frac{1}{2}\right],    
\end{equation}
\begin{equation}
W_2=\mbox{diag}\left[\frac{5}{12},\frac{13}{12},1,\ldots,1,\frac{13}{12},\frac{5}{12}\right],    
\end{equation}
\begin{equation}
W_3=\mbox{diag}\left[\frac{3}{8},\frac{7}{6},\frac{23}{24},1,\ldots,1,\frac{23}{24},\frac{7}{6},\frac{3}{8}\right],    
\end{equation}
\begin{equation}
W_4=\mbox{diag}\left[\frac{251}{720},\frac{299}{240},\frac{211}{240},\frac{739}{720},1,\ldots,1,\frac{739}{720},\frac{211}{240},\frac{299}{240},\frac{251}{720}\right]. 
\end{equation}
The index $n$ for the matrix $W_n$ means that $n$ weight coefficients are different from the unit. Gregory formulas $W_{2k-1}$ and $W_{2k}$ are accurate for polynomials up to $2k-1$ degree \cite{hamming1962numerical}. In our work, we use Gregory formulas up to $n=6$. One can use formulas with larger numbers. However, negative weight coefficients appear for $n=9$ and larger \cite{fornberg2019improved}.
If the solution decreases rapidly at infinity, then one-sided Gregory formulas can be used \cite{fornberg2019improved}.
}}

{\textcolor{black}{
As a result of the approximation for (\ref{X1P}), we obtain a system of linear equations with the Hankel matrices
\begin{equation}\label{Hb}
\begin{bmatrix}E&\mp h H^*W_n\\h H W_n&E\end{bmatrix}
\begin{bmatrix} \vec{X}_1\\ \vec{X}_2\end{bmatrix} =   
\begin{bmatrix} 0\\ \vec{F}\end{bmatrix},
\end{equation}
where $E$ is an identity matrix of size  $M+1$, vectors $\vec{X}_1$ and $\vec{X}_2$ have dimension $M+1$
$$\vec{X_1}=\left[X_1(0,t),X_1(h,t),X_1(2h,t),\ldots, X_1((M-1)h,t),X_1(Mh,t)\right]^T,$$ $$\vec{X_2}=\left[X_2(0,t),X_2(h,t),X_2(2h,t),\ldots,X_2((M-1)h,t), X_2(Mh,t)\right]^T,$$ $$\vec{F}=-\left[\omega_0,\omega_1,\ldots,\omega_M\right]^T.$$
}}
{\textcolor{black}{
Let's make a transformation $Y_m=W_n X_m$, $m=1,2$ and get
\begin{equation}\label{sysY}
\begin{bmatrix}W^{-1}_n&\mp h H^*\\h H&W_n^{-1}\end{bmatrix}
\begin{bmatrix}\vec{Y}_1\\ \vec{Y}_2\end{bmatrix}
=\begin{bmatrix}0\\ \vec{F}\end{bmatrix}.
\end{equation}
We will make another transformation to get a block system with Toeplitz matrices. Multiply system (\ref{sysY}) by block-diagonal matrix $\mbox{diag}(E,J)$
\begin{equation}\label{H2T}
\begin{bmatrix}E&0\\0&J\end{bmatrix}
\begin{bmatrix}W_n^{-1}&\mp h H^*\\h H&W_n^{-1}\end{bmatrix}
\begin{bmatrix}E&0\\0&J\end{bmatrix}
\begin{bmatrix}E&0\\0&J\end{bmatrix}
\begin{bmatrix}\vec{Y}_1\\\vec{Y}_2\end{bmatrix}=
\begin{bmatrix}E&0\\0&J\end{bmatrix}\begin{bmatrix}0\\\vec{F}\end{bmatrix},
\end{equation}
where $J$ is an exchange matrix \cite{blahut2010fast}. It is a matrix with the field element one in every entry
of the antidiagonal (the entries where $j = n+1-i$) and with the field element zero
in every other matrix entry. Notice that $J=J^T$, $J^2=E$ and the multiplication of the vector by the exchange matrix $J$ rearranges the elements of the vector in the reverse order. Simplifying expression (\ref{H2T}), we get a system with the matrix $B$
\begin{equation}\label{sysL}
\begin{bmatrix}W_n^{-1}&\mp T^*\\T&W_n^{-1}\end{bmatrix}
\begin{bmatrix}\vec{Y}_1\\J\vec{Y}_2\end{bmatrix}
=\begin{bmatrix}0\\J\vec{F}\end{bmatrix},\quad B=\begin{bmatrix}W_n^{-1}&\mp T^*\\T&W_n^{-1}\end{bmatrix},
\end{equation}
where the matrix $T$ and its Hermitian adjoint $T^*$ are Toeplitz matrixes and are expressed by $H$ and $J$
\begin{equation}
T=hJ H,\quad T^*=h H^*J.
\end{equation}
We can enter the notation for $J\vec{Y}_2$ and $J\vec{F}$, but in fact they are the vectors with the reverse order of elements.
}}

{\textcolor{black}{
The matrix $W_n$ is an almost identical matrix, so the matrix $B$ of the system (\ref{sysL}) is the difference of the block Toeplitz matrix $A$ and the diagonal matrix $R$
\begin{equation}
B=A-R,\quad A=\begin{bmatrix}E&\mp T^*\\T&E\end{bmatrix},\quad R=\begin{bmatrix}E-W_n^{-1}&0\\0&E-W_n^{-1}\end{bmatrix}.
\end{equation}
The rang $r$ of $R$ is equal to $4n$. In our situation $r<<M$, therefore we can apply the Woodbury formula \cite{hager1989updating}
\begin{equation}\label{woodbury}
\left(A-UV\right)^{-1}=A^{-1}+{A^{-1}U}\left(E_r-V\,{A^{-1}U}\right)^{-1}V\,A^{-1},
\end{equation}
where $A$, $E_r$, $U$, $V$ are conformable matrices and $r$ is a rang of $UV$,
for solving the linear system $B\vec{x}=\vec{b}$.
The matrix $R$ is diagonal therefore it is easy to factorize it by the general formula for the block matrix
\begin{equation}
\begin{bmatrix}
D_1&0&0&0&0&0\\0&0&0&0&0&0\\0&0&D_2&0&0&0\\
0&0&0&D_3&0&0\\0&0&0&0&0&0\\0&0&0&0&0&D_4
\end{bmatrix}=\begin{bmatrix}
U_1&0&0&0\\0&0&0&0\\0&U_2&0&0\\0&0&U_3&0\\
0&0&0&0\\0&0&0&U_4
\end{bmatrix}\begin{bmatrix}
V_1&0&0&0&0&0\\0&0&V_2&0&0&0\\
0&0&0&V_3&0&0\\0&0&0&0&0&V_4
\end{bmatrix},  
\end{equation}
where $D_k=U_kV_k$, $D_k$ is a non-singular diagonal matrix and zeros denotes conformable zero matrices. The simplest way to define $U_k$ and $V_k$ is to put $U_k=E_k$, $V_k=D_k$ or $U_k=D_k$, $V_k=E_k$. 
}}

{\textcolor{black}{
We now present a standard algorithm for solving a linear system $B\vec{x}=\vec{b}$. The solution $\vec{x}$ can be written out in the following form
\begin{equation}\label{woodbury2}
\vec{x}=\left(A-UV\right)^{-1}\vec{b}
=\underline{A^{-1}\vec{b}}+\underline{\underline{A^{-1}U}}
\left(E_r-V\underline{\underline{A^{-1}U}}\right)^{-1}V\,
\underline{A^{-1}\vec{b}}
\end{equation}
and finding it requires a sequence of 4 steps \cite{hager1989updating}:
\begin{enumerate}
    \item Find $\vec{y}$ by solving $A\vec{y}=\vec{b}$ (the underlined terms in (\ref{woodbury2})).
    \item Find the matrix $Z=A^{-1}U$ by solving $AZ=U$ (the twice underlined terms in (\ref{woodbury2})).
    \item Find $\vec{z}$ by solving $\left(E_r-VZ\right)\vec{z}=V\vec{y}$.
    \item Finally, $\vec{x}=\vec{y}+Z\vec{z}$.
\end{enumerate}
Thus, it is necessary to solve $r+1$ linear systems with the matrix $A$ of large size $M+1$ and one system with the matrix $\left(E_r-VZ\right)$ of small size $r$. The matrix $A$ has the Toeplitz blocks therefore one can apply a fast algorithm for the block Toeplitz matrix to solve this system.
}}

{\textcolor{black}{Let us describe the algorithm for the first step (see more details in \cite{MedvedevGTIB}). Let rewrite  the matrixes $T$, $T^*$ from $A$, the unknown vector $\vec{y}$ and the right hand size $\vec{b}$ 
\begin{equation}\label{matrix_form_T}
\begin{array}{l}
T^* = \begin{bmatrix}
		\omega_0^*   & \omega_1^* &  \ldots & \omega_{M}^*\\
		\omega_{-1}^* & \omega_0^* & &\vdots\\
		\vdots & &\rotatebox{00}{$\ddots$} &\vdots\\
		\omega_{-M}^* &\ldots &\ldots  & \omega_0^*
	\end{bmatrix},\quad
T =  \begin{bmatrix}
	\omega_0   & \omega_{-1} &  \ldots & \omega_{-M}\\
	\omega_{1} & \omega_0 & &\vdots\\
	\vdots & &\ddots &\vdots\\
	\omega_{M} & \ldots& \ldots & \omega_0
\end{bmatrix},\\[1cm]
\vec{y} = [y_0, y_1, \ldots, y_{M}, z_{0},z_1,\ldots, z_{M}]^T,\quad \vec{b} = [b_0, b_1,\ldots, b_{M}, c_{0},c_1\ldots, c_{M}]^T.
\end{array}
\end{equation}
The matrix of the system $A\vec{y}=\vec{b}$ consists of the Toeplitz blocks of size $M+1$. 
To solve this system numerically, we converted it to the system with the block-Toeplitz matrix with blocks of size $2$  \cite{ToeplitzGuide}: 
\begin{equation}\label{T2b:n}
\begin{bmatrix}
		t_0    & t_1 &  \ldots & t_{M}\\
		t_{-1} & t_0 & &\vdots\\
		\vdots & &\ddots &\vdots\\
		t_{-M} & \ldots& \ldots & t_0
\end{bmatrix}
\left[\begin{array}{c} x_0\\ \vdots \\ x_M\end{array}\right] =   
\left[\begin{array}{c} f_0\\ \vdots \\ f_M\end{array}\right],
\end{equation}
where
\begin{equation}
\begin{array}{l}
t_0 = \begin{bmatrix}
		1    & \mp \omega_0^*\\
		\omega_0 & 1 
	\end{bmatrix},\quad	
t_k = \begin{bmatrix}
		0    & \mp  \omega_{-k}^*\\
		\omega_k & 0 
	\end{bmatrix},\quad
t_{-k} = \begin{bmatrix}
		0    & \mp  \omega_k^*\\
		 \omega_{-k} & 0 
	\end{bmatrix},\quad
x_k = \left[\begin{array}{c} y_k\\ z_k\end{array}\right],\quad	
f_k = \left[\begin{array}{c} b_k\\ c_k\end{array}\right].
\end{array}
\end{equation}
To solve the system (\ref{T2b:n}) we used the block version of Levinson's algorithm \cite{ToeplitzGuide}. 
}}

{\textcolor{black}{This algorithm allows us to solve the system at a fixed point $t$ and, similarly to the procedure described in \cite{belai2007efficient}, allows us to find solutions of the GLME at the nearest point $t+h/2$. Since our approach is high order approximation and we can start calculations at an arbitrary point $t$ and then find solutions to the right of it with minimal computational costs, we named our method High-Order Generalized Toeplitz Inner-Bordering (HGTIB) method. 
{Similar} calculations for the right GLME allow finding solutions to the left of the starting point $t$ \cite{MedvedevGTIB}.
}}

\section{Numerical experiments}

We recover a potential $q(t)$ defined on a uniform grid of the interval of length $L$ with a step size $\tau = L/M$. 
Unless otherwise stated, the potential is recovered by solving the left GLME for $t \leq 0$ and the right GLME for $t > 0$.
Accordingly, the left GLME are solved on the interval $[-L/2, 0]$ and the right GLME are solved on the interval $(0, L/2]$. Each part corresponds to the integration region of size $P = L$ with a step size $h = 2\tau$.

Potential recovery error with respect to $t$ and the root-mean-square error are calculated in the entire computational interval by the formulas:
\begin{equation}\label{err}
\mbox{error}[q(t)] = \epsilon(t)\!=\!\frac{|q(t) - q^{exact}(t)|}{\max|q^{exact}(t)|},\quad
\mbox{RMSE}[q(t)]=\sqrt{\frac{\sum_{j=0}^{M}\epsilon(t_j)^2}{M+1}}
\end{equation}
\begin{figure}[htp]
\centering
\includegraphics[trim=1cm 2.5cm 5cm 0.5cm,clip=true, width=0.8\textwidth, draft=false]{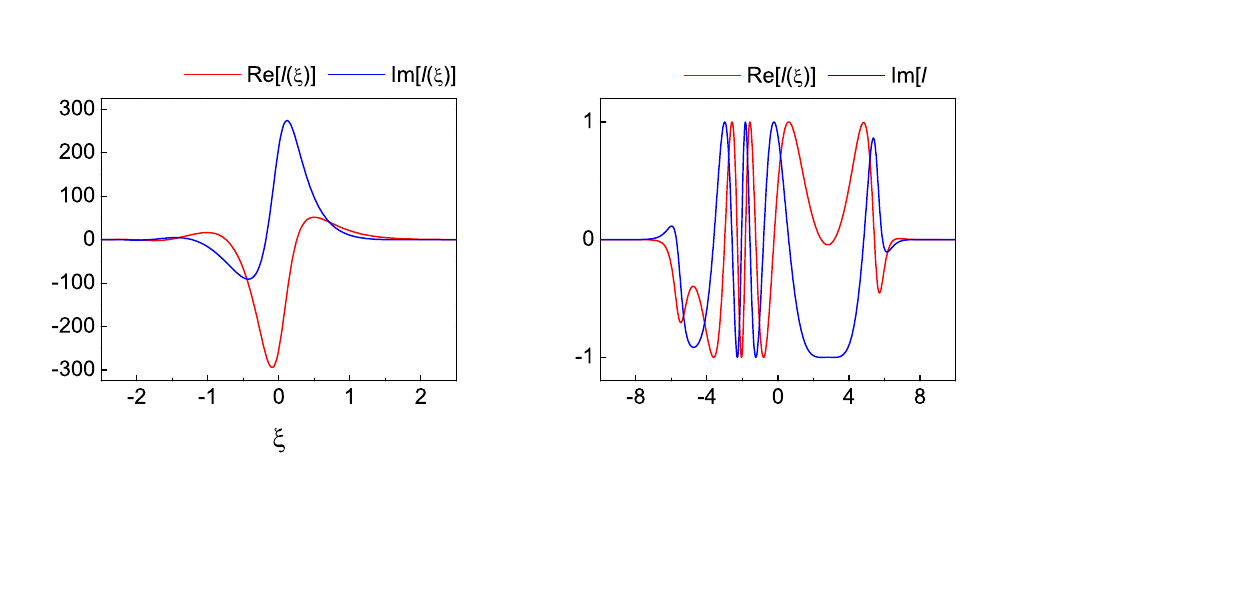}
\caption{Left reflection coefficient $l(\xi)$ for the chirped hyperbolic secant ($A = 5.2$, $C = 4$) in the case of anomalous (left) and normal (right) dispersion.}
\label{FigContSp}
\end{figure}
\begin{figure}[htp]
\centering
\includegraphics[trim=0.6cm 2.2cm 2.5cm 0.6cm,clip=true, width=0.9\textwidth, draft=false]{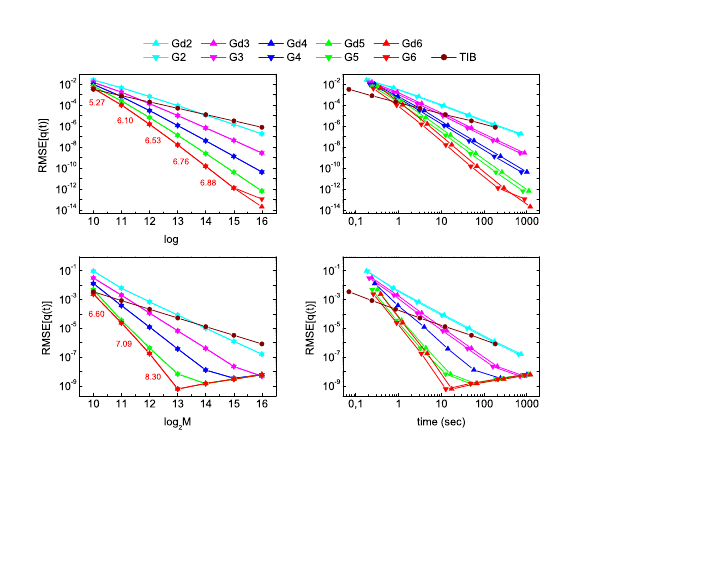}
\caption{Comparison of schemes with weight coefficients applying for both edges of the interval of integration
and schemes with weight coefficients applying only for one edge.
Root mean squared error (\ref{err}) with respect to the number of subintervals $M$ (left) and to the execution time trade-off (right)
in the case of anomalous (top row) and normal (bottom row) dispersion.}
\label{FigErr_vs_M}
\end{figure}
\begin{figure}[htp]
\centering
\includegraphics[trim=1cm 2.5cm 5cm 0.5cm,clip=true, width=0.8\textwidth, draft=false]{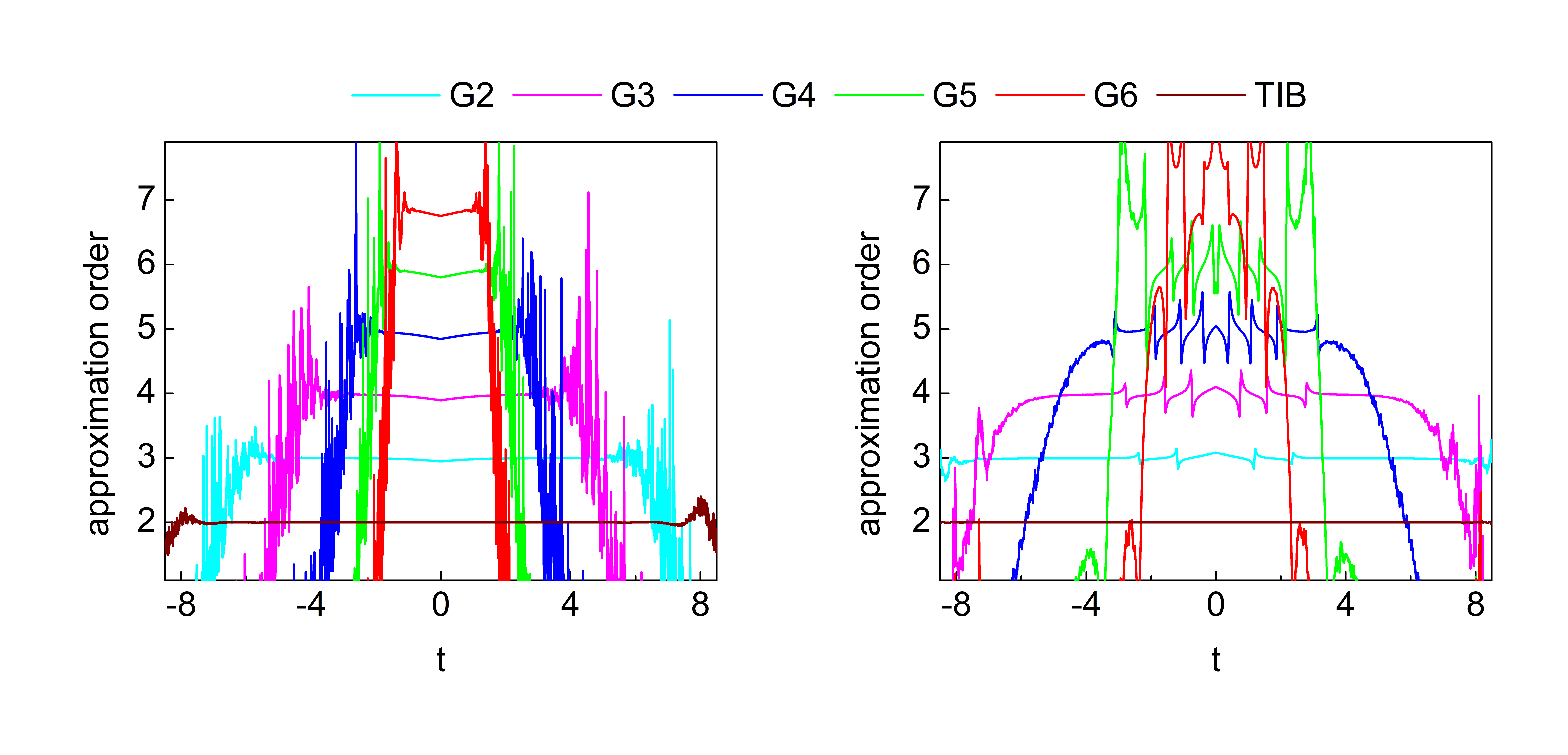}
\caption{Approximation order in the case of anomalous (left) and normal (right) dispersion.}
\label{FigApproxOrder}
\end{figure}
\begin{figure}[htp]
\centering
\includegraphics[trim=0.6cm 2.2cm 2.5cm 0.6cm,clip=true, width=0.9\textwidth, draft=false]{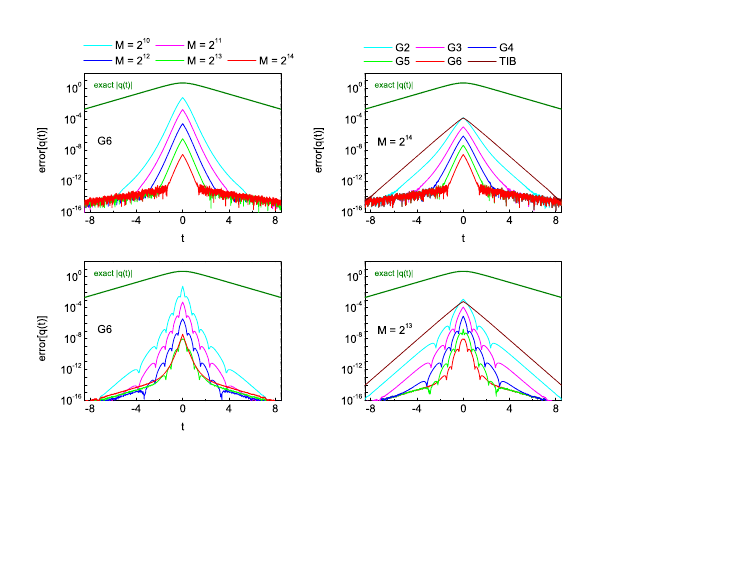}
\caption{The error $\epsilon(t)$ (\ref{err}) for different number of subintervals $M$ and different schemes in the case of anomalous (top row) and normal (bottom row) dispersion.}
\label{FigErr_vs_t}
\end{figure}


Here we compare the new high order schemes with the TIB scheme of the second order of accuracy.
Schemes G2d, G3d, G4d, G5d, G6d use approximation by the Gregory formula with $n = 2,3,4,5,6$ weight coefficients correspondingly.
Label "d" (double) means that weight coefficients are used in both edges of the interval of integration. In this case we need to invert a matrix of a rang $r =4n$.
We also present schemes G2, G3, G4, G5, G6 that use Gregory formula with applying weight coefficients only in one edge of the interval of integration. This means that we need to invert a matrix of a rang $r =2n$.

The potential in the form of a chirped hyperbolic secant $q(t) = A [\mbox{sech} (t)]^{1+iC}$ for $A = 5.2$, $C = 4$ is recovered by the combination of the left and right GLM equations. In the case of anomalous dispersion, such a signal has a continuous spectrum and five discrete eigenvalues. In the case of normal dispersion, it only has a continuous spectrum. The detailed 
analytical expressions of the spectral data for this type of potentials can be found in \cite{Medvedev2020_OE}.
To find the continuous part of the kernel (\ref{Omega_l}) we need to compute the integral for the reflection coefficient. 
We apply the trapezoidal rule if using the second order TIB scheme. When using higher order schemes to solve GLME, we apply the same Gregory's formula with wight coefficients always in both  edges of the interval of integration. In every numerical example here we use the same number of nodes $M_\xi = 2049$ in the spectral domain of size $L_\xi = 40$. Figure \ref{FigContSp} presents the continuous spectrum for this potential in the case of anomalous and normal dispersion.

Figure \ref{FigErr_vs_M} compares the schemes G2d, G3d, G4d, G5d, G6d with weight coefficients applying for both edges of the interval of integration
and schemes G2, G3, G4, G5, G6 with weight coefficients applying only for one edge.
Figure \ref{FigErr_vs_M} demonstrates the root mean squared error (\ref{err}) with respect to the number of subintervals $M$ (left) and 
to the execution time trade-off (right) in the case of anomalous and normal dispersion. 
One can see, that we can remove the weight coefficients from one edge of the interval of integration and not lose accuracy,
while the speed of the algorithm increases.
The best accuracy was provided by the schemes G6d and G6. 
Red numbers in Fig. \ref{FigErr_vs_M} (left) show the approximation order when moving from a grid with $M$ subintervals to a grid with $2M$ ones.
These numbers were calculated by the formula:
\begin{equation}\label{ApproxOrder}
m\!=\! \log_2\frac{\mbox{RMSE}[q(t)]^{2M}}{\mbox{RMSE}[q(t)]^M}.
\end{equation}
The arithmetic mean of the approximation order for G6 and G6d schemes is $6.31$ for the anomalous dispersion and $7.33$ for the normal dispersion. 

Figure \ref{FigErr_vs_M} (right) shows that the second order TIB scheme is the most efficient on coarse grids, but when one need to get the accuracy better than $10^{-4}$ the G6 scheme is the fastest.

Approximation order of the schemes with respect to $t$ is demonstrated in Fig. \ref{FigApproxOrder}. The order of approximation is calculated by a formula similar to (\ref{ApproxOrder}), but instead of the root-mean-square error, the error at each point $t$ is used. Here the number of subintervals $M = 2^{13}$ for anomalous and $M = 2^{12}$ for normal dispersion.

Figure \ref{FigErr_vs_t} (left) demonstrates the error $\epsilon(t)$ (\ref{err}) of recovering the potential $q(t)$ using the scheme G6 for different number of subintervals $M$ in the case of anomalous and normal dispersion. 
Figure \ref{FigErr_vs_t} (right) shows the same error $\epsilon(t)$ for different schemes G2 -- G6 and fixed number of subintervals $M$ for the anomalous and normal dispersion.

\section{Conclusion}
\label{s4}

A new high precision algorithm for solving the  Gelfand-Levitan-Marchenko equation is proposed. The algorithm is based on the block version of the Toeplitz Inner-Bordering algorithm of Livenson's type. To approximate integrals, we use the high-precision one-sided and two-sided Gregory quadrature formulas. Numerical experiments have shown that the use of the one-sided Gregory quadrature formulas does not lead to a loss of accuracy, but allows one to increase the speed of calculations.
The best accuracy was provided by the schemes G6d and G6, which use the Gregory formula with $6$ weight coefficients.
These schemes allows one to get the sixth approximation order for the anomalous dispersion and the seventh approximation order for the normal dispersion. Numerical experiments have also shown that the second order TIB scheme is the most efficient on coarse grids, but when one need to get the accuracy better than $10^{-4}$ the G6 scheme is the fastest. 

{\textcolor{black}{The proposed method can serve as part of an effective hybrid method in which the Gelfand-Levitan-Marchenko equation is used to find part of the potential for a continuous spectrum with high accuracy, and then the Darboux method is used to restore part of the potential for a discrete spectrum \cite{aref2018modulation}. }}

{\textcolor{black}{The proposed algorithm can be generalized to a vector version of the nonlinear Schrodinger equation. For example, for a two-component equation NLSE (Manakov equation), the Gelfand-Levitan-Marchenko equation has 3 components and its solution will be reduced to solving a system with a block Toeplitz matrix in which the blocks have 3 orders.
}}



\bibliography{main} 

\bibliographystyle{unsrt}

\end{document}